\documentclass{amsart}

\usepackage{amsmath,amsthm,amssymb}
\usepackage{amsfonts}
\usepackage[mathscr]{eucal}
\usepackage{ascmac}
\usepackage{epic, eepic}
\usepackage{graphicx}
\usepackage{cases}
\usepackage{color}
\usepackage{bm}
\usepackage[top=25truemm,bottom=25truemm,left=30truemm,right=30truemm]{geometry}
\usepackage{enumerate}

\theoremstyle{plane}
\newtheorem{thm}{Theorem}[section]
\newtheorem{prop}[thm]{Proposition}
\newtheorem{lem}[thm]{Lemma}
\newtheorem{cor}[thm]{Corollary}
\newtheorem{fact}[thm]{Fact}
\theoremstyle{definition}
\newtheorem{dfn}[thm]{Definition}

\theoremstyle{remark}
\newtheorem{rmk}[thm]{Remark}
\allowdisplaybreaks

\newcommand{\sgn}{\operatorname{sgn}}
\newcommand{\grad}{\operatorname{grad}}
\newcommand{\R}{\bm{R}}

\newcommand{\F}{\mathcal{F}}
\newcommand{\Sig}{\Sigma}
\numberwithin{equation}{section}
\renewcommand{\phi}{\varphi}
\renewcommand{\epsilon}{\varepsilon}
\newcommand{\what}{\widehat}
\newcommand{\wtilde}{\widetilde}

\begin{document}
%%%%% To ease editing, for IMPAN journals add:

%\baselineskip=17pt

%%%%%%%%%%%
\title[Focals of fronts]{Focal surfaces of wave fronts in the Euclidean 3-space}
\author[K. Teramoto]{Keisuke Teramoto}
\address{Department of Mathematics, Graduate School of Science, Kobe University, 
Rokko 1-1, Nada, Kobe 657-8501, Japan}
\email{teramoto@math.kobe-u.ac.jp}
%\date{July 26, 2016}

\maketitle

\begin{abstract}
We characterize singularities of focal surfaces of wave fronts 
in terms of differential geometric properties of the initial wave fronts. 
Moreover, we study relationships between geometric properties of focal surfaces 
and geometric invariants of the initial wave fronts.
\end{abstract}

\section{Introduction}
There is a class of surfaces with singularities called (wave) fronts. 
Recently, several differential geometric properties of fronts and criteria for singularities have been obtained 
(cf. \cite{fsuy,hhnsuy,im,krsuy,kruy,ms,msuy,mu,suy,suy2}). 
On the other hand, focal surfaces of regular surfaces have singularities in general. 
It is known that generic corank one singularities of focal surfaces are 
cuspidal edges and swallowtails (see \cite{agv,ifrt}).
Porteous \cite{po} found ridge points on regular surfaces, corresponding to cuspidal edges on focal surfaces.
Moreover, Bruce and Wilkinson \cite{bw} showed a relation between a sub-parabolic point 
and behavior of the Gaussian curvature of the focal surface (see also \cite{bgt1,bgt,morris}). 
For these reasons, studying focal surfaces may provide us with new geometric properties of surfaces. 
Thus we expect that we may obtain new geometric properties of fronts by investigating focal surfaces. 

In this paper, we study singularities and geometric properties of focal surfaces 
of fronts in the Euclidean $3$-space $\R^3$. 
It is known that one principal curvature is bounded and the other is unbounded 
at non-degenerate singular points (\cite{t2}, see also \cite{mu}). 
We show relations between singularities of the focal surface with respect to the bounded principal curvature 
and geometric properties of the original front (Theorem \ref{sing_caus}).

On the other hand, we also consider geometric properties of {the focal surface 
corresponding to the unbounded principal curvature.} 
This contains the image of the set of singular points of the original front as a curve on it. 
If the original front has a cuspidal edge, then the focal surface is regular near that cuspidal edge 
(Proposition \ref{regularity}). 
Thus we can consider the Gaussian curvature and the mean curvature of the focal surface corresponding to 
the cuspidal edge. 
We give explicit representations of the Gaussian and the mean curvature of the focal surface along the cuspidal 
edge by using geometric invariants (Theorem \ref{KH-hatc}). 

\section{Basic concepts and preliminaries}
In this section, we review some notions and concepts. 
For detailed explanations, see \cite{agv,ifrt,krsuy,msuy,suy3,suy,t1,t2}.
\subsection{Wave fronts}  
Let $\Sig$ be a domain of $(\R^2;u,v)$ and $f:\Sig\to\R^3$ a $C^\infty$ map. 
Then we call $f$ a {\it $($wave$)$ front} if there exists a unit vector field $\nu:\Sig\to S^2$ along $f$ such that 
\begin{itemize}
\item $\langle df_p(X_p),\nu(p)\rangle=0$ for any $p\in\Sig$, $X_p\in T_p\Sig$, and 
\item $L_f=(f,\nu):\Sig\to\R^3\times S^2$ gives an immersion, 
\end{itemize}
where $S^2$ is the unit sphere in $\R^3$ and 
$\langle\cdot,\cdot\rangle$ is the canonical inner product of $\R^3$. 
If $f$ satisfies the first condition, it is called a {\it frontal}. 
The vector field $\nu$ is called the {\it unit normal vector} or the {\it Gauss map} of $f$ (cf. \cite{agv,ifrt,krsuy}). 
A point $p\in\Sig$ is said to be a {\it singular point} if $f$ is not an immersion at $p$, and we denote by $S(f)$ 
the set of singular points of $f$.
 
For a frontal $f$, the function $\lambda:\Sig\to\R$ given by 
$$\lambda(u,v)=\det(f_u,f_v,\nu)(u,v)\quad (f_u=\partial f/\partial u,\ f_v=\partial f/\partial v)$$
is called the {\it signed area density function}. 
Then $S(f)=\{p\in\Sig|\lambda(p)=0\}$ holds. 
A singular point $p\in S(f)$ is said to be a {\it non-degenerate singular point} 
if $d\lambda(p)\neq0$. 
For a non-degenerate singular point $p$, there exist a neighborhood $V$ of $p$ and 
a regular curve $\gamma:(-\epsilon,\epsilon)\ni t\mapsto\gamma(t)\in V$ 
$(\epsilon>0)$ such that $\gamma(0)=p$ and $S(f)\cap V$ is parametrized by $\gamma$. 
Moreover, there exists a non-zero vector field $\eta$ on $S(f)\cap V$ 
such that $df_q(\eta_q)=0$ for $q\in S(f)\cap V$. 
We call $\gamma$ and $\eta$ the {\it singular curve} and the {\it null vector field}, respectively. 
The image $\hat{\gamma}=f\circ\gamma$ of a singular curve $\gamma$ via $f$ is called the {\it singular locus}. 
We denote by $\wtilde{\eta}$ the {\it extended null vector field} on $V$, 
which is a vector field defined on $V$ such that $df(\wtilde{\eta})=0$ on $S(f)\cap V$, i.e., 
$\wtilde{\eta}=\eta$ holds on $S(f)\cap V$ (see \cite{suy3}).

A non-degenerate singular point $p=\gamma(0)$ is said to be of the {\it first kind} if $\gamma'$ and $\eta$ 
are transverse at $t=0$. 
Otherwise, it is said to be of the {\it second kind} (\cite{msuy}). 
Moreover, a non-degenerate singular point of the second kind $p$ is called {\it admissible} 
if $p=\gamma(0)$ is of the second kind and $\det(\gamma',\eta)(t)\neq0$ for $t\neq0$. 
We call a singular point $p$ of a front an {\it admissible non-degenerate singular point} 
if $p$ is of the first kind or is of the admissible second kind. 
We note that the singular locus of a non-degenerate singular point of the first kind is a regular curve, 
but for the second kind, it has a singularity. 
\begin{dfn}
Let $f,g:(\R^2,0)\to(\R^3,0)$ be $C^\infty$ map germs. 
Then $f$ and $g$ are {\it $\mathcal{A}$-equivalent} if there exist diffeomorphism germs 
$\theta:(\R^2,0)\to(\R^2,0)$ and $\Theta:(\R^3,0)\to(\R^3,0)$ such that 
$\Theta\circ f=g\circ\theta$ holds.
\end{dfn}
\begin{dfn}
Let $f:(\Sig,p)\to(\R^3,f(p))$ be a $C^\infty$ map germ around $p$. 
Then $f$ at $p$ is a \textit{cuspidal edge} if the map germ $f$ is $\mathcal{A}$-equavalent to the map germ 
$(u,v)\mapsto(u,v^2,v^3)$ at ${0}$, and 
$f$ at $p$ is a \textit{swallowtail} if the map germ $f$ is $\mathcal{A}$-equivalent to the map germ 
$(u,v)\mapsto(u,3v^4+uv^2,4v^3+2uv)$ at ${0}$. 
%$f$ at $p$ is a \textit{cuspidal butterfly} if the map germ $f$ is $\mathcal{A}$-equivalent to the map germ 
%$(u,v)\mapsto(u,4v^5+uv^2,5v^4+2uv)$ at ${0}$  
\end{dfn}
A cuspidal edge is non-degenerate singular point of the first kind, and 
a swallowtail is of the admissible second kind. 
Criteria for these singularities are known.
\begin{fact}[{\cite{krsuy,suy3}}]
Let $f:\Sig\to\R^3$ be a front and $p$ a non-degenerate singular point. 
Then 
\begin{itemize}
\item[(1)] $f$ is a cuspidal edge at $p$ if and only if $\eta\lambda(p)\neq0$.
\item[(2)] $f$ is a swallowtail at $p$ if and only if $\eta\lambda(p)=0$ and $\eta\eta\lambda(p)\neq0$.
%\item[(3)] $f$ is a cuspidal butterfly at $p$ if and only if 
%$\eta\lambda(p)=\eta\eta\lambda(p)=0$ and $\eta\eta\eta\lambda(p)\neq0$.
\end{itemize}
\end{fact}
We note that criteria for other singularities are known as well (cf. \cite{is,ist,s}). 

\subsection{Geometric invariants of wave fronts}
We recall relations between geometric invariants and principal curvatures of fronts. 

Let $f:\Sig(\subset\R^2)\to\R^3$ be a front and $p\in\Sig$ a non-degenerate singular point. 
Then we can take the following local coordinate system around $p$.  
\begin{dfn}[\cite{krsuy,suy,msuy}]
Let $f:\Sig\to\R^3$ be a front and $p\in\Sig$ a non-degenerate singular point of the first kind 
(resp. the second kind). 
Then a local coordinate system $(U;u,v)$ centered at $p$ is an {\it adapted coordinate system} 
if the following properties hold:
\begin{enumerate}
\item the singular curve is the $u$-axis,
\item $\eta=\partial_v$ (resp. $\eta=\partial_u+e(u)\partial_v$ with $e(0)=0$) gives a null vector field
on the $u$-axis,
\item there are no singular points other than the $u$-axis.
\end{enumerate}
\end{dfn}

\begin{figure}[htbp]
  \begin{center}
    \begin{tabular}{c}

      % 1
      \begin{minipage}{0.33\hsize}
        \begin{center}
          \includegraphics[width=4.5cm]{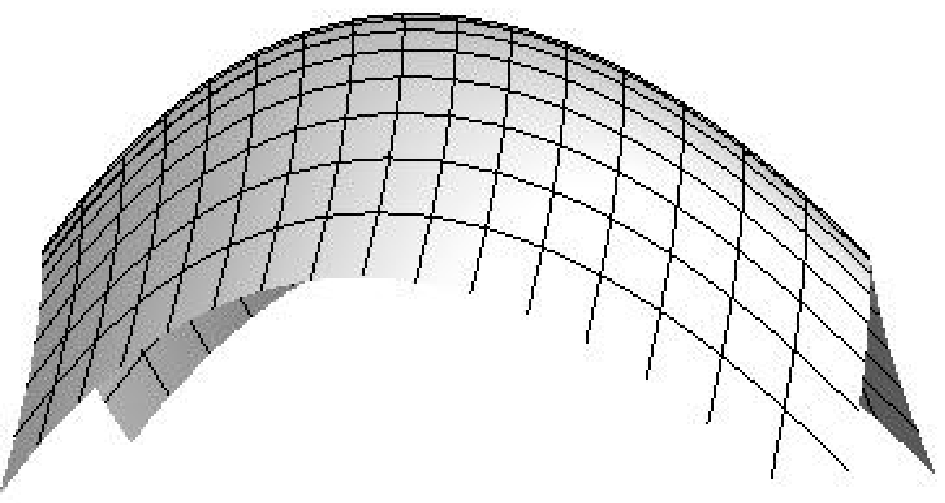}
         % \hspace{1cm} 
        \end{center}
      \end{minipage}

      % 2
      \begin{minipage}{0.33\hsize}
        \begin{center}
          \includegraphics[width=4.5cm]{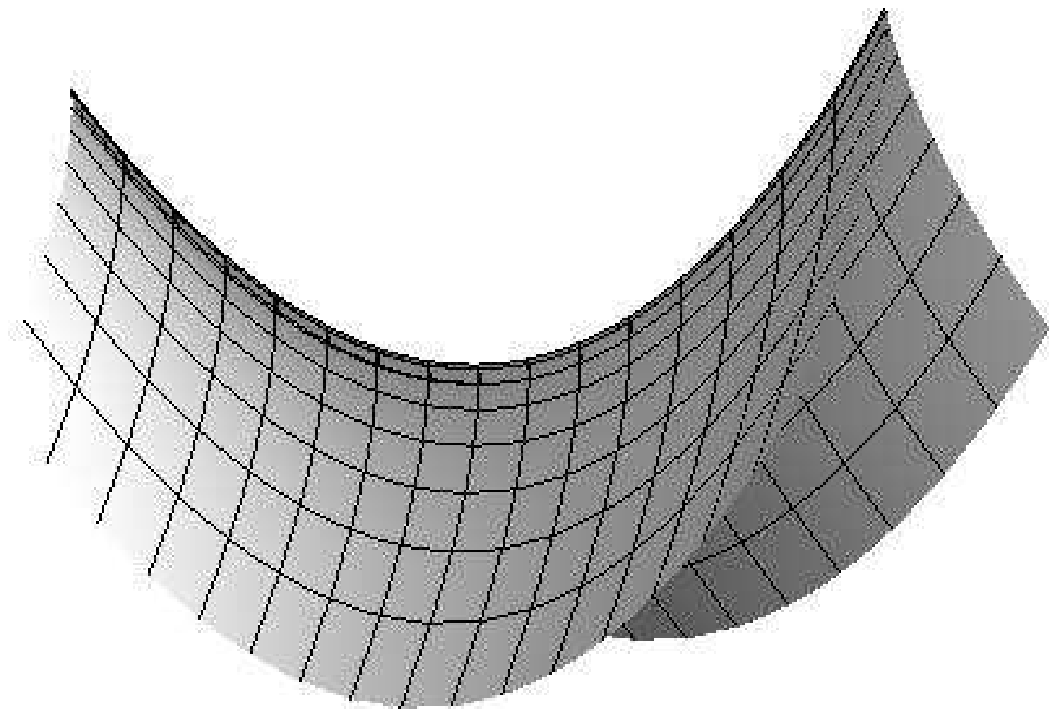}
          %\hspace{1cm} 
        \end{center}
      \end{minipage}

    \end{tabular}
    \caption{Cuspidal edges with positive $\kappa_s$ (left) and negative $\kappa_s$ (right).}
    \label{fig:kappas}
  \end{center}
\end{figure}

First we deal with cuspidal edges. 
Let $f:\Sig\to\R^3$ be a front, $\nu$ its unit normal vector and 
$p\in\Sig$ a cuspidal edge. 
Let $\kappa_s$, $\kappa_\nu$, $\kappa_c$ and $\kappa_t$ denote the {\it singular curvature} (\cite{suy}), 
the {\it limiting normal curvature} (\cite{suy}), the {\it cuspidal curvature} (\cite{msuy}) 
and the {\it cusp-directional torsion} (\cite{ms}), respectively. 
{If we take an adapted coordinate system $(U;u,v)$ around $p$, then 
\begin{align}\label{kappas}
\begin{aligned}
\kappa_s&=\sgn(\lambda_v)\dfrac{\det(f_u,f_{uu},\nu)}{|f_u|^3},\quad
\kappa_\nu=\dfrac{\langle{f_{uu},\nu}\rangle}{|f_u|^2},\quad
\kappa_c=\dfrac{|f_u|^{3/2}\det(f_u,f_{vv},f_{vvv})}{|f_u\times f_{vv}|^{5/2}},\\
\kappa_t&=\dfrac{\det(f_u,f_{vv},f_{uvv})}{|f_u\times f_{vv}|^2}
-\dfrac{\det(f_u,f_{vv},f_{uu})\langle{f_u,f_{vv}}\rangle}{|f_u|^2|f_u\times f_{vv}|^2}
\end{aligned}
\end{align}
hold on the $u$-axis, where $|\cdot|$ denotes the norm of $\R^3$.} 
We note that $\kappa_c$ does not vanish along the singular curve 
if it consists of cuspidal edges (\cite[Proposition 3.11]{msuy}). 
{Moreover, the sign of $\kappa_s$ corresponds to convexity 
or concavity of a cuspidal edge (see Figure \ref{fig:kappas}) (cf. \cite{suy}).}
For other geometric meanings of these invariants, see \cite{suy,msuy,ms}. 
Take an adapted coordinate system $(U;u,v)$ centered at $p$. 
Since $df(\eta)=f_v=0$ on the $u$-axis, there exists a $C^\infty$ map 
$h:U\to\R^3\setminus\{0\}$ such that $f_v=vh$ on $U$. 
We note that $\{f_u,h,\nu\}$ gives a frame on $U$, 
and we may take $\nu$ as $\nu=(f_u\times h)/|f_u\times h|$ (cf. \cite{ms,msuy,nuy}). 
Under the adapted coordinate system $(U;u,v)$ centered at $p$ with $\lambda_v(u,0)=\det(f_u,h,\nu)(u,0)>0$, 
$\kappa_\nu$, $\kappa_c$ and $\kappa_t$ can be written as follows (\cite[Lemma 2.7]{t2}): 
\begin{equation}\label{inv_cusp}
\kappa_\nu(u)=\frac{\wtilde{L}}{\wtilde{E}}(u,0),\quad
\kappa_c(u)=\frac{2\wtilde{E}^{3/4}\wtilde{N}}{(\wtilde{E}\wtilde{G}-\wtilde{F}^2)^{3/4}}(u,0),\quad 
\kappa_t(u)=\frac{\wtilde{E}\wtilde{M}-\wtilde{F}\wtilde{L}}{\wtilde{E}\sqrt{\wtilde{E}\wtilde{G}-\wtilde{F}^2}}(u,0),
\end{equation}
where 
\begin{equation}\label{fundamental}
\wtilde{E}=|f_u|^2,\quad \wtilde{F}=\langle f_u,h \rangle,\quad \wtilde{G}=|h|^2,\quad
\wtilde{L}=-\langle f_u,\nu_u\rangle,\quad \wtilde{M}=-\langle h,\nu_u\rangle,\quad\wtilde{N}=-\langle h,\nu_v\rangle.
\end{equation}
%Moreover, $\kappa_s$ is given by 
%\begin{equation}\label{kappas}
%\kappa_s(u)=\frac{\det(f_u,f_{uu},\nu)}{|f_u|^3}(u,0)
%\end{equation}
%if $\eta\lambda(u,0)>0$ holds (cf. \cite{suy,hhnsuy,ms}). 

We consider principal curvatures of cuspidal edges. 
We set two functions $\kappa_j:U\setminus\{v=0\}\to\R$ $(j=1,2)$ by 
\begin{equation}\label{principal}
\kappa_1=H+\sqrt{H^2-K},\quad \kappa_2=H-\sqrt{H^2-K},
\end{equation}
where $K$ and $H$ are the Gaussian and mean curvature of $f$ defined on 
$U\setminus\{v=0\}$, respectively.
Since $K=\kappa_1\kappa_2$ and $2H=\kappa_1+\kappa_2$ on $U\setminus\{v=0\}$, 
$\kappa_j$ $(j=1,2)$ can be considered as {\it principal curvatures} of $f$. 
Although $H$ is unbounded near $p$, it is known that $\kappa_2$ is bounded at $p$ 
if $\eta\lambda(p)\cdot\kappa_c(p)>0$ holds (\cite[Theorem 3.1]{t2}). 
Using functions defined in \eqref{fundamental}, $\kappa_1$ and $\kappa_2$ as in \eqref{principal} can be written as
\begin{equation}
\kappa_1=\frac{2(\wtilde{L}\wtilde{N}-v\wtilde{M}^2)}{\wtilde{A}-\wtilde{B}},\quad
\kappa_2=\frac{2(\wtilde{L}\wtilde{N}-v\wtilde{M}^2)}{\wtilde{A}+\wtilde{B}},\label{k1}
\end{equation}
where 
$\wtilde{A}=\wtilde{E}\wtilde{N}-2v\wtilde{F}\wtilde{M}+v\wtilde{G}\wtilde{L},\,
\wtilde{B}=\sqrt{\wtilde{A}^2-4v(\wtilde{E}\wtilde{G}-\wtilde{F}^2)(\wtilde{L}\wtilde{N}-v\wtilde{M}^2)}$. 

Suppose that $\kappa_2$ is bounded on $U$ and write $\kappa=\kappa_2$. 
Then $\kappa=\kappa_\nu$ holds along the $u$-axis (cf. \cite{t1,t2}). 
Moreover, we assume that $\kappa\neq0$ on $U$ in what follows, when we treat cuspidal edges. 
Let $\bm{V}=V_1\partial_u+V_2\partial_v$ be a vector field on $U$, 
where $(V_1,V_2)=(\wtilde{N}-v\kappa\wtilde{E},-\wtilde{M}+\kappa\wtilde{F})$. 
This vector field $\bm{V}$ gives a {\it principal vector} with respect to $\kappa$ (see \cite{t1,t2}). 

Next we treat a front with non-degenerate singular points of the second kind. 
Let $f:\Sig\to\R^3$ be a front, $\nu$ a unit normal vector to $f$ and 
$p\in \Sig$ a non-degenerate singular point of the second kind. 
Take an adapted coordinate system $(U;u,v)$  around $p$. 
Since $df(\eta)=0$ on the $u$-axis, there exists a $C^\infty$ map $h:U\to\R^3\setminus\{{0}\}$ such that 
$df(\eta)=f_u+e(u)f_v=vh$ on $U$. 
Since $d\lambda=\det(h,f_v,\nu)dv\neq0$ on the $u$-axis, $h,f_v$ and $\nu$ are linearly independent, 
and we may take an adapted coordinate system $(U;u,v)$ satisfying $\lambda_v(u,0)>0$. 
We define the following functions as 
\begin{equation}\label{fundamental2}
\what{E}=| h|^2,\quad\what{F}=\langle h,f_v \rangle,\quad\what{G}=|f_v|^2,\quad
\what{L}=-\langle h,\nu_u\rangle,\quad\what{M}=-\langle h,\nu_v\rangle,\quad\what{N}=-\langle f_v,\nu_v\rangle.
\end{equation}
We note that $\what{L}+e(u)\what{M}\neq0$ on the $u$-axis, 
in particular $\what{L}\neq0$ at $p$ (cf. \cite{t2}).

Let $K$ and $H$ denote the Gaussian and mean curvature of $f$ defined on $U\setminus\{v=0\}$. 
We set $\what{H}=vH$. 
It is known that $\what{H}$ is a $C^\infty$ function on $U$. 
Moreover, we set $\mu_c(p)=2\what{H}(p)$. 
This is a geometric invariant called the {\it normalized cuspidal curvature} (see \cite{msuy}). 
By using the functions above, we see that 
$\mu_c(p)=\what{G}(p)\what{L}(p)/|h(p)\times f_v(p)|^2$ holds (\cite{msuy,t2}). 
In addition, it is known that $\mu_c(p)\neq0$ if and only if $f$ is a front at $p$ (\cite[Proposition 3.2]{msuy}). 

We define functions $\kappa_j:U\setminus\{v=0\}\to\R$ $(j=1,2)$ as 
\begin{equation}\label{prinicipal2}
\kappa_1=H+\sqrt{H^2-K},\quad \kappa_2=H-\sqrt{H^2-K}.
\end{equation}
These functions can be considered as {\it principal curvatures} of $f$ on $U\setminus\{v=0\}$. 
It is known that if $\mu_c(p)>0$, then $\kappa_2$ is bounded at $p$ (\cite[Theorem 3.1]{t2}). 
Moreover, we can write $\kappa_1$ and $\kappa_2$ as 
\begin{equation}
\kappa_1=\frac{2((\what{L}+e(u)\what{M})\what{N}-v\what{M}^2)}{\what{A}-\what{B}},\quad
\kappa_2=\frac{2((\what{L}+e(u)\what{M})\what{N}-v\what{M}^2)}{\what{A}+\what{B}}\label{prinsec}
\end{equation}
by using functions in \eqref{fundamental2}, where 
\begin{equation*}
\what{A}=\what{G}(\what{L}+e(u)\what{M})-2v\what{F}\what{M}+v\what{E}\what{N},\quad%\label{sec_A}\\
\what{B}=\sqrt{\what{A}^2-4v(\what{E}\what{G}-\what{F}^2)
\left((\what{L}+e(u)\what{M})\what{N}-v\what{M}^2\right)}\label{sec_B}.
\end{equation*}
%$\kappa_j$ $(j=1,2)$ are also well-defined $C^\infty$ functions on $U\setminus\{v=0\}$. 
We assume that $\kappa_2$ is bounded  on $U$ and set $\kappa=\kappa_2$. 
Let $\bm{V}=V_1\partial_u+V_2\partial_v$ be a vector field on $U$ defined by 
$(V_1,V_2)=(-\what{M}+\kappa\what{F}, \what{L}-\kappa(v\what{E}-e(u)\what{F}))$. 
This gives a {\it principal vector} with respect to $\kappa$ (see \cite{t2}). 

We note that Murata and Umehara \cite{mu} introduced the notion called the 
{\it principal curvature maps} on fronts, which is the generalization of the notion of principal curvature functions. 
{The relationship between principal curvature maps and 
principal curvatures as in \eqref{k1} and \eqref{prinsec} is mentioned in \cite{t2}.} 
\begin{rmk}\label{rmk:hatkappa}
If $\kappa_2$ is bounded at a non-degenerate singular point $q$, 
then $\kappa_1$ is unbounded near $q$.
We set $\hat{\kappa}=\lambda\kappa_1$, where $\lambda=\det(f_u,f_v,\nu)$. 
If $p$ is a cuspidal edge, then 
$\hat{\kappa}={\wtilde{E}\wtilde{N}}/{|f_u\times h|}>0$ 
holds along the $u$-axis. 
Moreover, if $p$ is of the admissible second kind, 
then $\hat{\kappa}=\what{G}(\what{L}+e(u)\what{M})/|h\times f_v|>0$ on the $u$-axis. 
Thus $\hat{\kappa}$ is a $C^\infty$ function on $U$. 
\end{rmk}
\begin{dfn}\label{ridge}
Under the above settings, 
a point $p$ is called a \textit{ridge point} if $\bm{V}\kappa(p)=0$ holds, 
where $\bm{V}\kappa$ denotes the directional derivative of $\kappa$ with respect to $\bm{V}$. 
Moreover, a point $p$ is called a \textit{k-th order ridge point} if 
$\bm{V}^{(m)}\kappa(p)=0\ (1\leq m\leq k)$ and $\bm{V}^{(k+1)}\kappa(p)\neq0$ hold, 
where $\bm{V}^{(m)}\kappa$ means the $m$-th directional derivative of $\kappa$ with respect to $\bm{V}$. 
We call the set of ridge points the {\it ridge line}. 
\end{dfn}
Porteous \cite{po1} defined ridge points to investigate singularities of focal surfaces of regular surfaces. 
For more detailed properties of ridge points on regular surfaces, see \cite{bgt1,bgt,fh,ifrt,po1,po}. 
%We note that relations between ridge points and umbilic points on surfaces are studied in \cite{bgt}. 
%Moreover, ridge points of Whitney umbrellas (or cross caps) are defined in \cite{fh2}.  
%call the image of ridges the {\it ridge line}. 
%We assume that ridges are regular curve near $p$, 
%that is, $d(\bm{V}\kappa)(p)$ does not vanish. 

\subsection{Sub-parabolic points on cuspidal edges}
We now consider sub-parabolic points with respect to a bounded principal curvature at cuspidal edges. 
Let $f:\Sig\to\R^3$ be a front, $\nu$ its unit normal vector and $p\in\Sig$ a cuspidal edge. 
Take an adapted coordinate system $(U;u,v)$ centered at $p$ satisfying $\eta\lambda(u,0)>0$.  
Then we assume that $\kappa_2$ as in \eqref{k1} is bounded on $U$ and write $\kappa=\kappa_2$. 

Let us denote $\wtilde{\bm{V}}=(\wtilde{V}_1,\wtilde{V}_2)$. 
If $\wtilde{\bm{V}}$ is a principal vector with respect to $\kappa_1$ on $U\setminus\{v=0\}$, 
then $\wtilde{\bm{V}}$ satisfies the relation $(II-\kappa_1 I)\wtilde{\bm{V}}=0$, 
where 
$$I=\begin{pmatrix} \langle{f_u,f_u}\rangle & \langle{f_u,f_v}\rangle \\
\langle{f_u,f_v}\rangle & \langle{f_v,f_v}\rangle\end{pmatrix},\quad
II=\begin{pmatrix} -\langle{f_u,\nu_u}\rangle & -\langle{f_u,\nu_v}\rangle\\
-\langle{f_v,\nu_u}\rangle & -\langle{f_v,\nu_v}\rangle \end{pmatrix}.$$
By using functions as in \eqref{fundamental}, we can write 
\begin{equation}\label{pdirk-}
\begin{pmatrix}
\wtilde{L}-\kappa_1\wtilde{E} & v(\wtilde{M}-\kappa_1\wtilde{F})\\
v(\wtilde{M}-\kappa_1\wtilde{F}) & v(\wtilde{N}-v\kappa_1\wtilde{G})
\end{pmatrix}
\begin{pmatrix}
\wtilde{V}_1 \\ \wtilde{V}_2
\end{pmatrix}=
\begin{pmatrix}
0 \\ 0
\end{pmatrix}.
\end{equation}
%where $\wtilde{E}$, $\wtilde{F}$, $\wtilde{G}$, $\wtilde{L}$, $\wtilde{M}$ and $\wtilde{N}$ are given as in \eqref{fundamental}. 
This relation \eqref{pdirk-} is equivalent to the following equations, once multiplying by $\lambda$:
$$
\begin{cases}
(\lambda\wtilde{L}-\hat{\kappa}\wtilde{E})\wtilde{V}_1+v(\lambda\wtilde{M}-\hat{\kappa}\wtilde{F})\wtilde{V}_2=0 \\
(\lambda\wtilde{M}-\hat{\kappa}\wtilde{F})\wtilde{V}_1+(\lambda\wtilde{N}-v\hat{\kappa}\wtilde{G})\wtilde{V}_2=0,
\end{cases}
$$
where $\hat{\kappa}=\lambda\kappa_1$. 
Thus we may take $\wtilde{\bm{V}}=(\wtilde{V}_1,\wtilde{V}_2)$ as 
\begin{equation}\label{tilv1}
\wtilde{\bm{V}}=(v(\lambda\wtilde{M}-\hat{\kappa}\wtilde{F}),-\lambda\wtilde{L}+\hat{\kappa}\wtilde{E})
\quad \text{or} \quad
\wtilde{\bm{V}}=(\lambda\wtilde{N}-v\hat{\kappa}\wtilde{G},-\lambda\wtilde{M}+\hat{\kappa}\wtilde{F}).
\end{equation}
We note that one of above vectors in \eqref{tilv1} is well-defined on $U$, since we can write $\wtilde{\bm{V}}$ as 
$\wtilde{\bm{V}}=(0,\hat{\kappa}\wtilde{E})$ or $\wtilde{\bm{V}}=(0,\hat{\kappa}\wtilde{F})$ on the $u$-axis. 
In particular, $\wtilde{\bm{V}}=(0,\hat{\kappa}\wtilde{E})\neq(0,0)$ on the $u$-axis. 
Thus we can take the ({\it extended}\/) {\it principal vector} $\wtilde{\bm{V}}$ with respect to $\kappa_1$ as 
\begin{equation}\label{tilv}
\wtilde{\bm{V}}=(v(\lambda\wtilde{M}-\hat{\kappa}\wtilde{F}),-\lambda\wtilde{L}+\hat{\kappa}\wtilde{E})
\end{equation}
on $U$. 
We now identify the vector $\wtilde{\bm{V}}=(\wtilde{V}_1,\wtilde{V}_2)$ with the vector field 
$\wtilde{\bm{V}}=\wtilde{V}_1\partial_u+\wtilde{V}_2\partial_v$ on $U$. 
By \eqref{tilv}, we note that $\wtilde{\bm{V}}$ is parallel to the null vector field $\eta=\partial_v$ along the $u$-axis. 
\begin{dfn}
Under the above settings, a point $p$ is called the {\it sub-parabolic point} of $f$ if 
$\wtilde{\bm{V}}\kappa(p)=0$ holds, where $\wtilde{\bm{V}}\kappa$ means the directional derivative of $\kappa$
in the direction $\wtilde{\bm{V}}$.
\end{dfn}
For geometric meanings of sub-parabolic points on a regular surface, see \cite{bgt1,bt,bw,ifrt,morris}. 
\begin{prop}\label{sub-parab}
Let $f:\Sig\to\R^3$ be a front and $p\in\Sig$ a cuspidal edge. 
Assume that $\kappa=\kappa_2$ or $\kappa=\kappa_1$ is a bounded principal curvature near $p$. 
Then a point $p$ is a sub-parabolic point if and only if $4\kappa_t^2+\kappa_s\kappa_c^2=0$ at $p$.
\end{prop}
\begin{proof}
Let $(U;u,v)$ be an adapted coordinate system centered at $p$ satisfying $\eta\lambda(u,0)>0$. 
In this case, we may take $\nu$ as $\nu=(f_u\times h)/|f_u\times h|$. 
We assume that $\kappa_2$ is a bounded principal curvature of $f$ on $U$ and set $\kappa=\kappa_2$, 
that is, $\wtilde{N}>0$ on the $u$-axis. 
The directional derivative $\wtilde{\bm{V}}\kappa$ of $\kappa$ with respect to $\wtilde{\bm{V}}$ as in \eqref{tilv} is 
$\wtilde{\bm{V}}\kappa=\wtilde{V}_2\partial \kappa/\partial v$ at $p$. 
By \eqref{k1}, $\partial \kappa/\partial v$ can be written as 
$$\frac{\partial \kappa}{\partial v}=2\left(\frac{\wtilde{L}_v\wtilde{N}+\wtilde{L}\wtilde{N}_v-
\wtilde{M}^2}{\wtilde{A}+\wtilde{B}}
-\frac{\wtilde{L}\wtilde{N}(\wtilde{A}_v+\wtilde{B}_v)}{(\wtilde{A}+\wtilde{B})^2}\right)$$
on the $u$-axis.
We note that $\wtilde{A}+\wtilde{B}=2\wtilde{E}\wtilde{N}$, 
$\wtilde{A}_v=\wtilde{E}\wtilde{N}_v-2\wtilde{F}\wtilde{M}+\wtilde{G}\wtilde{L}$ 
and $\wtilde{A}_v+\wtilde{B}_v=2(\wtilde{A}_v-\wtilde{L}|f_u\times h|^2/\wtilde{E})$ hold on the $u$-axis, 
since $2\wtilde{E}_v(u,0)=\langle f_{uv},f_u\rangle(u,0)=0$. 
Thus it follows that 
$$\frac{\partial \kappa}{\partial v}=\frac{\wtilde{L}_v\wtilde{N}-\kappa_t^2|f_u\times h|^2}{\wtilde{E}\wtilde{N}}$$
holds along the $u$-axis, by \eqref{inv_cusp}. 

We consider $\wtilde{L}_v$ along the $u$-axis. 
Since $\wtilde{L}=\det(f_u,h,f_{uu})/|f_u\times h|$ and $f_{uv}=f_{uuv}=0$ on the $u$-axis, we see that 
$$\wtilde{L}_v=\frac{\det(f_u,h_v,f_{uu})}{|f_u\times h|}-\frac{|f_u\times h|_v\det(f_u,h,f_{uu})}{|f_u\times h|^2}$$
holds at $p$. 
Since $h_v$ and $|f_u\times h|_v$ can be written as 
\begin{equation*}
h_v=\ast f_u+
\frac{|f_u|^2\langle h,h_v\rangle-\langle f_u,h\rangle\langle f_u,h_v\rangle}{|f_u\times h|^2}h+\wtilde{N}\nu,\quad
|f_u\times h|_v=\frac{|f_u|^2\langle h,h_v\rangle-\langle f_u,h\rangle\langle f_u,h_v\rangle}{|f_u\times h|},
\end{equation*}
we have 
$\wtilde{L}_v%={\what{N}\det(f_u,\nu,f_{uu})}{|f_u\times h|}
=-{\kappa_s\wtilde{N}|f_u|^3}/{|f_u\times h|}$ (see \eqref{kappas}). 
Thus $4\wtilde{L}_v\wtilde{N}= 
-{\kappa_s\kappa_c^2}|f_u\times h|^2$ 
holds on the $u$-axis, by \eqref{inv_cusp}. 
Hence we get 
$$\frac{\partial \kappa}{\partial v}
=-\frac{1}{2\kappa_c}(4\kappa_t^2+\kappa_s\kappa_c^2)\left(\frac{|f_u\times h|}{|f_u|}\right)^{1/2}$$
at $p$. 
Since $2\hat{\kappa}\wtilde{E}={\kappa_c}(|f_u|^5|f_u\times h|)^{1/2}$ at $p$, 
we have the assertion.
%$$\wtilde{\bm{V}}\kappa=\wtilde{\bm{V}}_2\kappa_v=-\frac{1}{4}(4\kappa_t^2+\kappa_s\kappa_c)(|f_u|^2|f_u\times h|)$$
%at $p$. 
\end{proof}

\section{Focal surfaces of wave fronts}\label{main}
We consider focal surfaces of fronts. 
Since focal surfaces can be regarded as the singular value set of a certain map $\R^3\to\R^3$, 
we consider focal surfaces by using results about Morin singularities of $\R^3\to\R^3$. 

\subsection{Morin singularities}
We recall relations between the $A_k$-Morin singularities and the $A_k$-front singularities.
The {\it $A_{k}$-Morin singularities} are map germs $f:(\R^n,p)\to(\R^n,f(p))$
which are $\mathcal{A}$-equivalent to 
$$f(x_1,\ldots,x_n)=(x_1,\ldots,x_{n-1},x_1x_n+\cdots+x_{k-1}x_n^{k-1}+x_n^{k+1})\quad (k\leq n)$$
at the origin $0$ (see \cite{morin,gg,suy3}). 
We note that the $A_0$-Morin singularity is actually a regular point.
\begin{fact}[{\cite[Theorem A.1]{suy3}}]\label{crit:morin}
Assume that $k\leq n$. 
Let $\Omega$ be a domain of $\R^n$, $f:\Omega\to\R^n$ a $C^\infty$ map and $p$ a singular point of $f$. 
Assume that $p$ is a corank one singularity. 
Then $f$ at $p$ is $\mathcal{A}$-equivalent to an $A_k$-Morin singularity if and only if 
\begin{itemize}
\item[(1)] $\Lambda=\Lambda'=\cdots=\Lambda^{(k-1)}=0$ and $\Lambda^{(k)}\neq0$ at $p$,
\item[(2)] $(\Lambda,\Lambda',\ldots,\Lambda^{(k-1)}):\Omega\to\R^k$ is non-singular at $p$,
\end{itemize}
where $\Lambda=\det(f_{x_1},\ldots,f_{x_n})$, $(x_1,\ldots,x_n)$ is the canonical coordinate system on $\Omega$, 
$\Lambda'=\wtilde{\eta}\Lambda$, $\Lambda^{(i)}=\wtilde{\eta}\Lambda^{(i-1)}$ and 
$\wtilde{\eta}$ is the extended null vector field of $f$.
\end{fact}
For stating the next fact, recall that an {\it $A_{k+1}$-front singularity} is a $C^\infty$ map germ defined as 
$$ X\mapsto\left((k+1)t^{k+2}+\sum_{j=2}^k(j-1)t^jx_j,-(k+2)t^{k+1}-\sum_{j=2}^kjt^{j-1}x_j,X_1\right)$$
at $0$, where $X=(t,x_2,\ldots,x_n)$ and $X_1=(x_2,\ldots,x_n)$ (see \cite{agv,suy3}). 
\begin{fact}[{\cite[Corollary 2.11]{suy3}}]\label{morin_front}
Let $\Omega$ be a domain of $\R^{n+1}$ and $f:\Omega\to\R^{n+1}$ a $C^\infty$ map. 
Suppose that $p\in\Omega$ is a singular point of $f$ such that the exterior derivative of the Jacobian of $f$ 
does not vanish at $p$. 
Then the following are equivalent:
\begin{itemize}
\item[(1)] $p$ is an $A_k$-Morin singular point of $f$,
\item[(2)] $f|_{S(f)}$ is a front, and $p$ is an $A_k$-front singularity of $f|_{S(f)}$.
\end{itemize}
\end{fact}
\begin{rmk}
The image of an $A_1$-front singularity is a regular point, 
the image of an $A_2$-front singularity is a cuspidal edge, 
and the image of an $A_3$-front singularity is a swallowtail 
if the dimensions of the source space and target space are two and three, respectively (see \cite{agv,bgt,suy3}).
\end{rmk}

\subsection{Focal surfaces of wave fronts}
Let $f:\Sig\to\R^3$ be a front, $\nu$ a unit normal vector to $f$ and $p\in \Sig$ a non-degenerate singular point. 
We assume that $p$ is of the second kind, and $(U;u,v)$ is an adapted coordinate system centered at $p$ 
satisfying $\lambda_v(u,0)>0$. 
When $p$ is a cuspidal edge, the following arguments can be applied similarly.

We now consider the map $\F:U\times\R\to\R^3$ given by 
\begin{equation}\label{bigf}
\F(u,v,w)=f(u,v)+w\nu(u,v) \quad\left((u,v)\in U,\ w\in\R\right).
\end{equation}
%This is a $C^\infty$ map on $U\times\R$. 
%It is known that $\F$ gives a {\it Lagrangian mapping} (\cite{agv,ifrt,istake}). 
%In fact, consider the map $\mathcal{L}=(\F,\nu):U\times\R\to T^\ast\R^3$, 
%where $T^\ast\R^3$ with local coordinates $(q_1,q_2,q_3,p_1,p_2,p_3)$ is the cotangent bundle over $\R^3$ 
%equipped with the canonical symplectic form 
%$\omega=\sum_{i=1}^3dp_i\wedge dq_i$. 
%Then $\dim (U\times\R)=3$, $\mathcal{L}$ gives an immersion and $\mathcal{L}^\ast\omega=0$ holds. 
%Thus $\mathcal{L}$ gives a {\it Lagrangian immersion}. 
%Let $\pi:T^\ast\R^3\to\R^3$ be a canonical projection. 
%Then it follows that $\F=\pi\circ\mathcal{L}$ holds, and hence $\F$ is a Lagrangian mapping. 
By direct computations, it follows that 
$$\F_u=(v+w\alpha_1)h+(-e(u)+w\alpha_2)f_v,\quad\F_v=w\beta_1h+(1+w\beta_2)f_v,\,\quad\F_w=\nu,$$
where 
\begin{equation*}
\alpha_1=\frac{\what{F}(v\what{M}-e(u)\what{N})-\what{G}\what{L}}{|h\times f_v|^2},\ 
\alpha_2=\frac{\what{F}\what{L}-\what{E}(v\what{M}-e(u)\what{N})}{|h\times f_v|^2},\ 
\beta_1=\frac{\what{F}\what{N}-\what{G}\what{M}}{|h\times f_v|^2},\ 
\beta_2=\frac{\what{F}\what{M}-\what{E}\what{N}}{|h\times f_v|^2}
\end{equation*}
(cf. \cite[Lemmas 2.6 and 2.8]{t2}). 
From these calculations, 
the Jacobian of $\F$ can be written as 
$$\det(\F_u,\F_v,\F_w)=(1-w\kappa)(\lambda-w\hat{\kappa}),$$
where $\hat{\kappa}=\lambda\kappa_1$ and $\lambda=\det(f_u,f_v,\nu)$. 
By Remark \ref{rmk:hatkappa}, 
$\hat{\kappa}$ is a $C^\infty$ function and does not vanish on the $u$-axis, in particular at $p$.
Thus we see that $\det(\F_u,\F_v,\F_w)(u,v,w)=0$ if and only if 
$1-w\kappa(u,v)=0$ or $\lambda-w\hat{\kappa}(u,v)=0$. 
Hence the set of singular points of $\F$ is $S(\F)=S_1(\F)\cup S_2(\F)$, where 
$S_1(\F)=\{(u,v,w)\,|\,w=1/\kappa(u,v)\}$ and $S_2(\F)=\{(u,v,w)\,|\,w=\lambda(u,v)/\hat{\kappa}(u,v)\}$. 
The image of $\F(S(\F))$ is
\begin{align}\label{img:s(f)}
\F(S(\F))&=
\left\{f(u,v)+\frac{1}{\kappa(u,v)}\nu(u,v)\,\Big|\,(u,v)\in U, w=\frac{1}{\kappa(u,v)}\right\}\\
&\bigcup\left\{f(u,v)+\frac{\lambda(u,v)}{\hat{\kappa}(u,v)}\nu(u,v)\,\Big|\,
(u,v)\in U, w=\frac{\lambda(u,v)}{\hat{\kappa}(u,v)}\right\}.\notag
\end{align}
We set 
\begin{equation}\label{caustics}
FC_f(u,v)=f(u,v)+\frac{1}{\kappa(u,v)}\nu(u,v),\quad \what{FC}_f(u,v)=f(u,v)+\frac{\lambda(u,v)}{\hat{\kappa}(u,v)}\nu(u,v).
\end{equation}
These are {\it focal surfaces} of $f$ (cf. \cite[pages 231 and 232]{doc}, see also \cite{cr,istake}). 
If $f$ at $p$ is a cuspidal edge, by similar calculations, we have the same formulae as in \eqref{caustics} 
for focal surfaces of $f$, 
where $\kappa=\kappa_2$ is given by \eqref{k1} and $\hat{\kappa}=\lambda\kappa_1$ 
is given by Remark \ref{rmk:hatkappa}. 
We assume that $f$ at $p$ is a cuspidal edge or of the second kind. 
We note that $FC_f$ cannot be defined at $p$ if $\kappa(p)=0$ (such a point is called a {\it parabolic point}). 
On the other hand, $\what{FC}_f$ can be defined near $p$ even if $\kappa(p)$ vanishes. 
Since the set of singular points $S(f^t)$ of a {\it parallel surface} $f^t=f+t\nu$, 
where $t\in\R\setminus\{0\}$ is constant, of a front $f$ is given by $S(f^t)=\{q\in U|t=1/\kappa(q)\}$,
the union of all the sets of singular points of the $f^t$ 
corresponds to the focal surface $FC_f$ if $\kappa$ never vanishes.  

We note that focal surfaces (or caustics) of regular surfaces in the Minkowski $3$-space $\R^{2,1}$ 
are studied in \cite{tari} by using distance squared functions on surfaces. 
{We also remark that Kokubu, Rossman, Umehara and Yamada \cite{kruy} studied 
caustics of flat fronts in the hyperbolic $3$-space $H^3$ from the global viewpoint.}

\subsection{Singularities of a focal surface $FC_f$ on a wave front}
We consider relations between singularities of $FC_f$ at $p$ and geometric properties of $f$. 
We assume that $p$ is not a parabolic point with respect to $\kappa$. 
\begin{lem}\label{corank}
Under the above settings, a singular point $P=(p,w_0=1/\kappa(p))\in S_1(\F)$ of $\F$ is corank one. 
Moreover, $S_1(\F)$ is a smooth submanifold of $U\times\R$ with codimension one near $P$.
\end{lem}
\begin{proof} 
We show the case that $p$ is a non-degenerate singular point of the second kind. 
For cuspidal edges, one can show this in a similar way. 

By the above calculations, $\F_w=\nu$ is linearly independent of $\F_u$ and $\F_v$. 
We note that $\F_u$ and $\F_v$ do not vanish at $P$ simultaneously, 
since $\alpha_1(p)\neq0$ holds. 
The cross product of $\F_u$ and $\F_v$ satisfies 
$$\F_u\times \F_v=(1-w\kappa)(\lambda-w\hat{\kappa})\nu=0$$
at $P=(p,w_0)\in S_1(\F)$. 
Thus $\F_u$ and $\F_v$ are linearly dependent at $P$. 
This implies that a point $P=(p,w_0)\in S_1(\F)$ is corank one. 

We show $S_1(\F)$ is a smooth submanifold of $U\times\R$ near $P$. 
By straightforward computations, the Jacobian matrix of $\F$ is rank two at $P=(p,w_0)\in S_1(\F)$. 
We set $\Lambda:U\times\R\to\R$ to be $\Lambda(u,v,w)=1-w\kappa(u,v)$. 
The gradient vector $\grad(\Lambda)$ of $\Lambda$ is
$$\grad(\Lambda)=\left(-\frac{\kappa_u}{\kappa},-\frac{\kappa_v}{\kappa},-\kappa\right)\neq(0,0,0)$$
at $(p,w_0)\in S_1(\F)$ since $\kappa(p)\neq0$, 
where $\kappa_u=\partial \kappa/\partial u$ and $\kappa_v=\partial \kappa/\partial v$. 
By the implicit function theorem, we have the conclusion.
\end{proof}

Let $\bm{V}$ be a principal vector with respect to $\kappa$. 
Then $d\F(\bm{V})=0$ holds on $S_1(\F)$, using the definitions of principal curvatures 
and principal vectors (cf. \cite{t1,t2}). 
Therefore $\bm{V}$ can be considered as the extended null vector field $\wtilde{\eta}$ of $\F$. 
\begin{lem}\label{lem:morin}
Under the above conditions, the following assertions hold.
\begin{itemize}
\item[(1)] $\F$ has an $A_1$-Morin singularity at $P=(p,w_0)\in S_1(\F)$ if and only if $p$ is not a ridge point of $f$.
\item[(2)] $\F$ has an $A_2$-Morin singularity at $P=(p,w_0)\in S_1(\F)$ if and only if 
$p$ is a first order ridge point of $f$.
\item[(3)] $\F$ has an $A_3$-Morin singularity at $P=(p,w_0)\in S_1(\F)$ if and only if 
$p$ is a second order ridge point of $f$ and the ridge line passing through $p$ is a regular curve.
\end{itemize}
Here $w_0=1/\kappa(p)$.
\end{lem}
\begin{proof}
%Let $f:U\to\R^3$ be a front, $p\in U$ a non-degenerate singular point of the second kind 
%and $\nu$ a unit normal vector to $f$. 
Let $\F:U\times\R\to\R^3$ be a $C^\infty$ map given by \eqref{bigf}. 
By Lemma \ref{corank}, it follows that a singular point $P=(p,w_0)\in S_1(\F)$ of $\F$ is corank one. 
Moreover, the extended null vector field $\wtilde{\eta}$ can be taken as 
a principal vector $\bm{V}$ of $\kappa$ for $\F$, 
and the function $\Lambda$ which gives $S_1(\F)$ can be taken as $\Lambda=1-w\kappa$ 
in both the case that $f$ at $p$ is a cuspidal edge and 
the case that it is a non-degenerate singular point of the second kind. 

First, we show $(1)$. 
Since $\kappa(p)\neq0$, $d\Lambda(P)\neq0$. 
By the assumptions, we have $\wtilde{\eta}\Lambda=-\bm{V}\kappa/\kappa\neq0$ at $p$ 
if and only if the point $p$ is not a ridge point. 
Thus assertion $(1)$ holds by Fact \ref{crit:morin}. 

Next, we prove the assertion $(2)$. 
We assume that $\wtilde{\eta}\Lambda=-\bm{V}\kappa/\kappa=0$ at $P=(p,w_0)$. 
The second order directional derivative of $\Lambda$ in the direction $\wtilde{\eta}$ becomes 
$\wtilde{\eta}^{(2)}\Lambda=-\bm{V}^{(2)}\kappa/\kappa$ 
%\wtilde{\eta^{(3)}}\Lambda=-\frac{\bm{V}^{(3)}\kappa_2}{\kappa_2}$$
at $P=(p,w_0)$. 
Moreover, the map $(\Lambda,\wtilde{\eta}\Lambda):U\times\R\to\R^2$ 
is non-singular at $P=(p,w_0)$ if and only if the matrix 
\begin{equation*}
\begin{pmatrix}
-\kappa_u/\kappa&-\kappa_v/\kappa&-\kappa\\
-(\bm{V}\kappa)_u/\kappa& -(\bm{V}\kappa)_v/\kappa&0
\end{pmatrix}
\end{equation*}
has rank two at $P=(p,w_0)$, by Fact \ref{crit:morin} in the case of $k=2$. 
Since $\bm{V}^{(2)}\kappa(p)\neq0$, $d(\bm{V}\kappa)(p)$ does not vanish, 
and hence assertion $(2)$ holds.

Finally, we show $(3)$. 
We assume that $\Lambda=\wtilde{\eta}\Lambda=\wtilde{\eta}^{(2)}\Lambda=0$ at $P=(p,w_0)\in S_1(\F)$, 
that is, $w_0=1/\kappa(p)$ and $\bm{V}\kappa(p)=\bm{V}^{(2)}\kappa(p)=0$ hold. 
Then $\wtilde{\eta}^{(3)}\Lambda\neq0$ at $P$ if and only if $\bm{V}^{(3)}\kappa(p)\neq0$. 
In addition, a map $(\Lambda,\wtilde{\eta}\Lambda,\wtilde{\eta}^{(2)}\Lambda):U\times\R\to\R^3$ 
is non-singular at $P$ if and only if the matrix
\begin{equation*}
\begin{pmatrix}
-\kappa_u/\kappa&-\kappa_v/\kappa&-\kappa\\
-(\bm{V}\kappa)_u/\kappa& -(\bm{V}\kappa)_v/\kappa&0\\
-(\bm{V}^{(2)}\kappa)_u/\kappa& -(\bm{V}^{(2)}\kappa)_v/\kappa&0
\end{pmatrix}
\end{equation*}
has rank three at $P$. 
Since $\bm{V}^{(3)}\kappa(p)\neq0$, $d(\bm{V}^{(2)}\kappa)$ does not vanish at $p$. 
Therefore the above $3\times3$ matrix has rank three at $P$ if and only if $d(\bm{V}\kappa)$ does not vanish at $p$. 
This condition is equivalent to the condition that the ridge line passing through $p$ is a regular curve. 
Thus we have the assertion by Fact \ref{crit:morin}. 
\end{proof}

For the focal surface $FC_f$, we shall prove the following assertion.
\begin{thm}\label{sing_caus}
Let $f:\Sig\to\R^3$ be a front and $p\in\Sig$ a non-degenerate singular point. 
Suppose that $\kappa=\kappa_2$ $($resp. $\kappa=\kappa_1)$ is a $C^\infty$ principal curvature of $f$ near $p$ 
and $FC_f$ is a focal surface of $f$ with respect to $\kappa$. 
Then the following assertions hold.
\begin{itemize}
\item[(1)] $FC_f$ is non-singular at $p$ if and only if $p$ is not a ridge point of $f$.
\item[(2)] $FC_f$ is a cuspidal edge at $p$ if and only if $p$ is a first order ridge point of $f$.
\item[(3)] $FC_f$ is a swallowtail at $p$ if and only if $p$ is a second order ridge point of $f$ and 
the ridge line passing through $p$ is a regular curve.
\end{itemize}
\end{thm}
\begin{proof}
We prove the case that the front has a non-degenerate singular point of the second kind. 
For the case of cuspidal edges, we can show this in a similar way. 

Let $f:\Sig\to\R^3$ be a front, $\nu$ a unit normal vector of $f$ 
and $p\in \Sig$ a non-degenerate singular point. 
Suppose that $(U;u,v)$ is an adapted coordinate system centered at $p$ and 
a principal curvature $\kappa$ of $f$ is of class $C^\infty$ on $U$. 
We construct a map $\F:U\times\R\to\R^3$ as in \eqref{bigf}. 
Then the image of the set of singular points of $\F$ gives a focal surface of $f$ 
with respect to $\kappa$ (see \eqref{img:s(f)}). 
Moreover, the point $P=(p,w_0)\in S_1(\F)$ is a corank one singular point of $\F$, by Lemma \ref{corank}. 
Thus we get the conclusions 
by Lemma \ref{lem:morin} and Fact \ref{morin_front}.
\end{proof}

\subsection{Geometric properties of $\what{FC}_f$ of cuspidal edges}\label{other}
We consider geometric properties of $\what{FC}_f$ as in \eqref{caustics} of a front $f:\Sig\to\R^3$ 
with a cuspidal edge $p\in\Sig$. 
For a front $f:\Sig\to\R^3$ with unit normal $\nu$ and a cuspidal edge point at $p\in\Sig$, 
the {\it limiting tangent plane of $f$ at $f(p)$} is the plane which is perpendicular to $\nu(p)$.
\begin{prop}\label{regularity}
Let $f:\Sig\to\R^3$ be a front, $\nu$ a unit normal vector to $f$ and $p\in\Sig$ a cuspidal edge. 
Then the focal surface $\what{FC}_f$ is regular at $p$. 
Moreover, the limiting tangent plane $LT$ of $f$ at $f(p)$ and the tangent plane $TP$ of $\what{FC}_f$ 
at $\what{FC}_f(p)$ 
intersect orthogonally.
\end{prop}

\begin{proof}
Let us take an adapted coordinate system $(U;u,v)$ around $p$ satisfying $\eta\lambda(u,0)>0$. 
Then the limiting tangent plane $LT$ of $f$ at $f(p)$ is generated by $f_u(p)$ and $h(p)$. 
Moreover, $\nu$ is given by $\nu=(f_u\times h)/|f_u\times h|$. 

On the other hand, we consider the tangent plane of the focal surface $\what{FC}_f$ at $\what{FC}_f(p)$. 
We note that $\what{FC}_f(p)=f(p)$ holds. 
By direct computations, we have 
\begin{equation*}
(\what{FC}_f)_u=f_u+\left(\frac{\lambda}{\hat{\kappa}}\right)_u\nu+\frac{\lambda}{\hat{\kappa}}\nu_u,\quad
(\what{FC}_f)_v=vh+\left(\frac{\lambda}{\hat{\kappa}}\right)_v\nu+\frac{\lambda}{\hat{\kappa}}\nu_v.
\end{equation*}
Thus $(\what{FC}_f)_u(p)=f_u(p)$ and $(\what{FC}_f)_v(p)=\lambda_v(p)\nu(p)/\hat{\kappa}(p)$ hold, 
where $\lambda_v=\det(f_u,h,\nu)=|f_u\times h|$. 
This implies that $\what{FC}_f$ is regular at $p$, and 
$f_u(p)$ and $\nu(p)$ are an orthogonal basis of the tangent plane $TP$ of $\what{FC}_f$ at $\what{FC}_f(p)$. 
A normal vector to $\what{FC}_f$ is given as $\wtilde{\bm{n}}=f_u\times\nu$ along the $u$-axis. 
Since $\langle\nu,\wtilde{\bm{n}}\rangle=0$ holds on the $u$-axis, 
$LT$ and $TP$ intersect orthogonally at $f(p)=\what{FC}_f(p)$.
\end{proof}
By Proposition \ref{regularity}, we can consider the Gaussian and mean curvature of $\what{FC}_f$ along the 
singular curve $\gamma$ of $f$. 
\begin{thm}\label{KH-hatc}
The Gaussian curvature $K_{\what{FC}_f}$ and mean curvature $H_{\what{FC}_f}$ 
of the focal surface $\what{FC}_f$ are given as 
$$K_{\what{FC}_f}=-\frac{1}{4}(4\kappa_t^2+\kappa_s\kappa_c^2),\quad 
H_{\what{FC}_f}=\pm\frac{1}{8}(\kappa_c^2-4\kappa_s)$$
along $\gamma$, 
where the sign $\pm$ of $H_{\what{FC}_f}$ depends on the orientation of the unit normal vector to $\what{FC}_f$.
\end{thm}
\begin{proof}
Let us take an adapted coordinate system $(U;u,v)$ centered at $p$ with $\eta\lambda(u,0)=\lambda_v(u,0)>0$. 
Then we may take $\nu$ as $\nu=(f_u\times h)/|f_u\times h|$, 
and we have $\lambda_v=\det(f_u,h,\nu)=|f_u\times h|$. 
Since $(\what{FC}_f)_u=f_u$ and $(\what{FC}_f)_v=\lambda_v\nu/\hat{\kappa}$ on the $u$-axis, 
the coefficients of the first fundamental form of $\what{FC}_f$ are 
$$E_{\what{FC}_f}=\wtilde{E}=|f_u|^2,\quad F_{\what{FC}_f}=0,\quad G_{\what{FC}_f}=\frac{\lambda_v^2}{\hat{\kappa}^2}$$
along $\gamma$. 
The second order differentials of $\what{FC}_f$ can be written as 
$$(\what{FC}_f)_{uu}=f_{uu},\quad(\what{FC}_f)_{uv}=\frac{\lambda_v}{\hat{\kappa}}\nu_u+\ast_1\nu,
\quad(\what{FC}_f)_{vv}=-h+\frac{\ast_2}{\hat{\kappa}}f_u+\ast_3\nu$$
on $\gamma$, where $h:U\to\R^3\setminus\{0\}$ is a $C^\infty$ map satisfying 
$f_v=vh$ and $\ast_i$ ($i=1,2,3$) are some functions. 
We can take a unit normal vector $\bm{n}$ to $\what{FC}_f$ as 
$\bm{n}=\pm(f_u\times\nu)/|f_u|$ along $\gamma$. 
Thus coefficients of the second fundamental form of $\what{FC}_f$ are 
$$L_{\what{FC}_f}=\pm\frac{\det(f_u,\nu,f_{uu})}{|f_u|},\quad
M_{\what{FC}_f}=\pm\frac{\lambda_v\det(f_u,\nu,\nu_{u})}{\hat{\kappa}|f_u|},\quad 
N_{\what{FC}_f}=\pm\frac{\det(f_u,h,\nu)}{|f_u|}$$
at $p$. 
By \eqref{kappas}, $L_{\what{FC}_f}=\mp\kappa_s|f_u|^2$ holds.
By \cite[Lemma 2.1]{t1} and \eqref{inv_cusp}, $\nu_u$ is expressed as 
\begin{equation}\label{nuu}
\nu_u=\frac{\wtilde{F}\wtilde{M}-\wtilde{G}\wtilde{L}}{|f_u\times h|^2}f_u-\frac{\kappa_t|f_u|^2}{|f_u\times h|}h
\end{equation}
along $\gamma$. 
On the other hand, the following equation holds on $\gamma$, by \eqref{inv_cusp}:
\begin{equation}\label{hatk-kc}
\frac{\lambda_v}{\hat{\kappa}}=\frac{|f_u\times h|^2}{|f_u|^2\wtilde{N}}=
\frac{2}{\kappa_c}\left(\frac{|f_u\times h|}{|f_u|}\right)^{1/2}.
\end{equation}
Hence $M_{\what{FC}_f}$ is calculated as 
$$M_{\what{FC}_f}=\pm\frac{2}{\kappa_c}\left(\frac{|f_u\times h|}{|f_u|}\right)^{1/2}
\frac{\kappa_t|f_u|}{|f_u\times h|}\det(f_u,h,\nu)
=\pm\frac{2\kappa_t|f_u|^{1/2}}{\kappa_c|f_u\times h|^{1/2}}\det(f_u,h,\nu)$$
on $\gamma$, by \eqref{nuu} and \eqref{hatk-kc}. 
Since $\det(f_u,h,\nu)=|f_u\times h|$, we have 
$$M_{\what{FC}_f}=\pm\frac{2\kappa_t}{\kappa_c}(|f_u||f_u\times h|)^{1/2},\quad
N_{\what{FC}_f}=\pm\frac{|f_u\times h|}{|f_u|}$$
and  
\begin{multline*}
E_{\what{FC}_f}G_{\what{FC}_f}-F_{\what{FC}_f}^2=\frac{4|f_u||f_u\times h|}{\kappa_c^2},\quad
L_{\what{FC}_f}N_{\what{FC}_f}-M_{\what{FC}_f}^2=
-\frac{(4\kappa_t^2+\kappa_s\kappa_c^2)}{\kappa_c^2}|f_u||f_u\times h|,\\ 
E_{\what{FC}_f}N_{\what{FC}_f}-2F_{\what{FC}_f}M_{\what{FC}_f}+G_{\what{FC}_f}L_{\what{FC}_f}=
\pm\frac{(\kappa_c^2-4\kappa_s)}{\kappa_c^2}|f_u||f_u\times h|
\end{multline*}
along $\gamma$. 
Thus the assertions hold by the following formulae:
$$K_{\what{FC}_f}=\frac{L_{\what{FC}_f}N_{\what{FC}_f}-M_{\what{FC}_f}^2}
{E_{\what{FC}_f}G_{\what{FC}_f}-F_{\what{FC}_f}^2}\quad\text{and}\quad
H_{\what{FC}_f}=\frac{E_{\what{FC}_f}N_{\what{FC}_f}-2F_{\what{FC}_f}M_{\what{FC}_f}+G_{\what{FC}_f}L_{\what{FC}_f}}
{E_{\what{FC}_f}G_{\what{FC}_f}-F_{\what{FC}_f}^2}.$$
\end{proof}
Comparing Theorem \ref{KH-hatc} and Proposition \ref{sub-parab}, we have the following assertion.
\begin{cor}\label{vanishing}
Let $f$ be a front in $\R^3$, $p\in\Sig$ a cuspidal edge of $f$ and $\what{FC}_f$ the focal surface. 
Then the Gaussian curvature $K_{\what{FC}_f}$ of $\what{FC}_f$ vanishes at $p$ 
if and only if $p$ is a sub-parabolic point with respect to a bounded principal curvature of $f$.
\end{cor}
This property is similar to the case of regular surfaces obtained by Buruce and Wilkinson \cite{bw} and 
Morris \cite{morris} (see also \cite{ifrt}).

For the focal surface $\what{FC}_f$, 
we call a point $p$ an {\it elliptic}, a {\it parabolic} or a {\it hyperbolic} point of $\what{FC}_f$ if 
$K_{\what{FC}_f}>0$, $=0$ or $<0$ at $p$, respectively.
By Theorem \ref{KH-hatc}, we have the following properties immediately.
\begin{cor}\label{characteristic}
Let $f$ be a front in $\R^3$, $p$ a cuspidal edge of $f$ and $\what{FC}_f$ a focal surface of $f$.  
Then we have the following.
\begin{enumerate}
\item[$(1)$] The point $p$ of $\what{FC}_f$ is classified as follows$:$
\begin{itemize}
\item $p$ is an elliptic point of $\what{FC}_f$ if and only if $4\kappa_t^2+\kappa_s\kappa_c^2<0$ at $p$,
\item $p$ is a parabolic point of $\what{FC}_f$ if and only if $4\kappa_t^2+\kappa_s\kappa_c^2=0$ at $p$,
\item $p$ is a hyperbolic point of $\what{FC}_f$ if and only if $4\kappa_t^2+\kappa_s\kappa_c^2>0$ at $p$.
\end{itemize}
Moreover, the Gaussian curvature $K_{\what{FC}_f}$ is non-negative 
at $p$ if and only if $\kappa_s$ is non-positive at $p$. 
In particular, if $K_{\what{FC}_f}$ is strictly positive along $\gamma$, 
then $\kappa_s$ is strictly negative. 
\item[$(2)$] If the mean curvature $H_{\what{FC}_f}$ vanishes, 
$\kappa_s$ is strictly positive along $\gamma$.
\end{enumerate}
\end{cor}
The invariant $4\kappa_t^2+\kappa_s\kappa_c^2$ appears as the coefficient of $v$ in 
the Gaussian curvature $K$ of a cuspidal edge (see \cite{msuy}).

Under this setting, since the singular locus $\hat{\gamma}$ is a regular curve on $\what{FC}_f$, 
we can consider the {\it geodesic curvature} $\hat{\kappa}_g$ and the {\it normal curvature} $\hat{\kappa}_n$ 
of $\what{FC}_f$ along the singular curve $\gamma$. 
\begin{prop}\label{kgkn}
Let $f:\Sig\to\R^3$ be a front, $p\in\Sig$ a cuspidal edge of $f$, $\gamma$ a singular curve passing through $p$ 
and $\eta$ a null vector field. 
Assume that $\eta\lambda>0$ $($resp. $\eta\lambda<0)$ along $\gamma$. 
Then $\hat{\kappa}_g=\kappa_\nu$ and $\hat{\kappa}_n=-\kappa_s$ 
$($resp. $\hat{\kappa}_g=\kappa_\nu$ and $\hat{\kappa}_n=\kappa_s)$ hold along $\gamma$. 
\end{prop}
\begin{proof}
Let us take an adapted coordinate system $(U;u,v)$ centered at $p$ with $\eta\lambda(u,0)>0$. 
Then we take the unit normal vector $\bm{n}$ to $\what{FC}_f$ along $\gamma$ as 
$\bm{n}=(f_u\times\nu)/|f_u|$. 
The geodesic curvature $\hat{\kappa}_g$ and the normal curvature $\hat{\kappa}_n$ of $\what{FC}_f$ along $\gamma$ 
are written as 
$$\hat{\kappa}_g=\frac{\langle \hat{\gamma}'',\bm{n}\times\hat{\gamma}'\rangle}{|\hat{\gamma}'|^3}
=\frac{\langle f_{uu},\bm{n}\times f_u\rangle}{|f_u|^3},\quad 
\hat{\kappa}_n=\frac{\langle\hat{\gamma}'',\bm{n}\rangle}{|\hat{\gamma}'|^2}=\frac{\langle f_{uu},\bm{n}\rangle}{|f_u|^2}.
$$
By direct calculations, we see that 
$\bm{n}\times f_u=|f_u|\nu$ 
and $\langle f_{uu},\bm{n}\rangle=-\det(f_u,f_{uu},\nu)/|f_u|$ hold. 
By \cite[(1.7)]{suy} and \cite[(3.11)]{suy}, we have the assertions.
\end{proof}
We recall the (pre-)geodesics on a surface. 
A curve on a regular surface is called a {\it pre-geodesic} if the geodesic 
curvature vanishes along the curve (cf. \cite{klin}). 
In addition, we call a curve on a regular surface a {\it geodesic} if the curve is pre-geodesic and has unit speed. 
\begin{cor}\label{pregeo}
The Gaussian curvature of a cuspidal edge is bunded on a sufficiently 
small neighborhood of the singular curve $\gamma$ if and only if 
$\hat{\gamma}$ is a $($pre-$)$geodesic on $\what{FC}_f$.
\end{cor}
\begin{proof}
The Gaussian curvature of a cuspidal edge is bounded if and only if 
the limiting normal curvature $\kappa_\nu$ vanishes along 
the singular curve $\gamma$ (\cite[Theorem 3.1]{suy}). 
Thus we have the assertion by Proposition \ref{kgkn}.
\end{proof}

It is known that the singular locus $\hat{\gamma}$ of cuspidal edges is a {\it line of curvature} 
if and only if the cusp-directional torsion $\kappa_t$ vanishes identically on $\gamma$ 
(\cite[Proposition 3.2]{t2}, see also \cite{istake2}).
In this case, we have the following.

\begin{prop}\label{curvline}
Let $f$ be a front, $p$ a cuspidal edge of $f$ and $\gamma$ a singular curve passing through $p$. 
Suppose that $\hat{\gamma}$ is a line of curvature on $f$. 
Then $\hat{\gamma}$ is also a line of curvature on $\what{FC}_f$.
\end{prop}

\begin{proof}
It is known that a curve $\sigma(t)$ on a regular surface is a line of curvature if and only if 
$$\det(\dot\sigma,n,\dot{n})=0$$
holds, where $n=n(t)$ is a unit normal to the surface restricted to $\sigma$ and 
we denote $\dot{~}=d/dt$. 
We apply this fact to the case of $\what{FC}_f$. 

We take an adapted coordinate system $(U;u,v)$ centered at $p$ satisfying $\eta\lambda(u,0)>0$. 
Then the unit normal vector to $\what{FC}_f$ can be taken as $\bm{n}=(f_u\times\nu)/|f_u|$ along $\gamma(u)=(u,0)$. 
Differentiating $\bm{n}$, we have 
$$\bm{n}'=\bm{n}_u=\frac{f_{uu}\times\nu+f_u\times\nu_u}{|f_u|}-|f_u|_u(f_u\times\nu).$$
Thus $\det(\hat{\gamma}',\bm{n},\bm{n}')$ can be written as 
\begin{equation*}
\det(\hat{\gamma}',\bm{n},\bm{n}')=\frac{1}{\wtilde{E}}\det(f_u,f_u\times\nu,f_{uu}\times\nu+f_u\times\nu_u)
=\frac{1}{\wtilde{E}}\det(f_u,f_u\times\nu,f_u\times\nu_u),
\end{equation*}
since $\langle f_u,\nu\rangle=0$, 
where we used the relation 
$\det(\bm{a}\times \bm{b},\bm{a}\times \bm{c},\bm{d})=\det(\bm{a},\bm{b},\bm{c})\langle{\bm{a},\bm{d}}\rangle$ 
($\bm{a},\bm{b},\bm{c},\bm{d}\in\R^3$). 
By \eqref{nuu}, $\nu_u$ is written as 
$$\nu_u=\frac{\wtilde{F}\wtilde{M}-\wtilde{G}\wtilde{L}}{|f_u\times h|^2}f_u-\frac{\kappa_t|f_u|^2}{|f_u\times h|}h.$$
Hence we have
\begin{equation*}
\det(\hat{\gamma}',\bm{n},\bm{n}')
%=\frac{\wtilde{F}\wtilde{L}-\wtilde{E}\wtilde{M}}{\wtilde{E}(\wtilde{E}\wtilde{G}-\wtilde{F}^2)}
%\det(f_u,f_u\times\nu,f_u\times h)=
=\frac{\kappa_t\wtilde{E}}{\sqrt{\wtilde{E}\wtilde{G}-\wtilde{F}^2}}\det(f_u,h,\nu).
\end{equation*} 
Since $\lambda_v=\det(f_u,h,\nu)$ does not vanish on $\gamma$, we have the assertion. 
\end{proof}
In general, a line of curvature on a regular surface does not become 
a line of curvature on its focal surface (cf. \cite[Proposition 6.19]{ifrt}). 
Thus Proposition \ref{curvline} gives a characteristic of cuspidal edges. 

\proof[Acknowledgements]
The author would like to express his sincere gratitude to Professor Kentaro Saji 
for his constant encouragement and fruitful discussions, 
and to Professor Wayne Rossman for his advice and helpful comments. 
The author also would like to thank the referee for reading this manuscript carefully and offering useful comments.
This work was supported by Grant-in-Aid for JSPS Research Fellows Number 17J02151.

%\vspace{0.5cm}
%\small{\begin{flushleft}
%\begin{tabular}{l}
%Keisuke Teramoto\\
%Department of Mathematics, \\
%Kobe University,\\
%Rokko 1-1, Nada, \\
%Kobe 657-8501, Japan\\
%E-mail adress: {\tt teramotoO\!\!\!amath.kobe-u.ac.jp}\\
%\end{tabular}
%\end{flushleft}}

\end{document}